\documentclass[11pt]{article}
\catcode`\@=11 \catcode`\@=12
\usepackage[utf8]{inputenc}
\usepackage[T1]{fontenc}
\usepackage{amsmath}
\usepackage{amssymb}
\usepackage{amsthm}
\usepackage{oldgerm}
\usepackage{multicol}
\usepackage{verbatim}
\usepackage{color}
\usepackage{ulem}
\usepackage{enumitem}

\usepackage[pdftex,colorlinks,urlcolor=blue,pdfstartview=FitH]{hyperref}

\makeatletter
\newcommand{\mylabel}[2]{#2\def\@currentlabel{#2}\label{#1}}
\makeatother

\newcommand{\Rm}{\mathbb{R}}

\newcommand{\RR}{\mathbb{R}}

\newcommand{\mC}{\ensuremath{\mathcal{C}}}
\newcommand{\mR}{\ensuremath{\mathcal{R}}}

\newcommand{\mF}{\ensuremath{\mathcal{F}}}

\newcommand{\Sm}{\ensuremath{\mathbb{S}}}

\newcommand{\Zm}{\ensuremath{\mathbb{Z}}}

\newcommand{\mK}{\ensuremath{\mathcal{K}}}
\newcommand{\mD}{\ensuremath{\mathcal{D}}}

\newcommand{\mI}{\ensuremath{\mathcal{I}}}
\newcommand{\mJ}{\ensuremath{\mathcal{J}}}

\newcommand{\mE}{\ensuremath{\mathcal{E}}}

\newcommand{\vs}{\vspace{.2cm}}

\newcommand{\gC}{\ensuremath{\bf C}}
\newcommand{\gL}{\ensuremath{\bf L}}
\newcommand{\gS}{\ensuremath{\bf S}}

\newtheorem{lem}{Lemma}

\newtheorem{thm}{Theorem}

\newtheorem{cor}[lem]{Corollary}
\newtheorem{prop}[lem]{Proposition}
\newtheorem{defn}[lem]{Definition}

\newtheorem{exmp}[lem]{Example}

\def\proof {\noindent{\sc{Proof. }}}
\def\qed {\mbox{}\hfill {\small \fbox{}} \\}
\def\lto{\longrightarrow}

\def\leq{\leqslant}
\def\geq{\geqslant}

\begin{document}

\begin{center}
	\begin{huge}
		{\bf Smoothing causal functions.}\\
	\end{huge}
	\vs
-----
\vs
\end{center}
\begin{small}
\begin{multicols}{2}

\noindent
Patrick Bernard
\footnote{Universit\'e  Paris-Dauphine},\\
PSL Research University,\\
\'Ecole Normale Sup\'erieure,\\
CNRS, DMA (UMR 8553)\\
45, rue d'Ulm\\
75230 Paris Cedex 05,
France\\
\texttt{patrick.bernard@ens.fr}\\

\noindent
Stefan Suhr,\\
Fakult\"at f\"ur Mathematik,\\
Ruhr-Universit\"at Bochum\\
Universit\"atsstra\ss e 150,\\
44801 Bochum,\\
Germany\\
\texttt{stefan.suhr@rub.de}\\

\end{multicols}
\vs
\thispagestyle{empty}
\begin{center}
-----
\end{center}

\textbf{Abstract. } 
We describe, in the general setting of closed cone fields, the set of causal functions
which can be approximated by smooth Lyapunov.
We derive several consequences on causality theory.

\begin{center}
-----
\end{center}

\textbf{R\'esum\'e. }
Dans le contexte général des champs de cones fermés, on décrit l'ensemble des fonctions
causales qui peuvent être approchées par des fonctions de Lyapunov lisses.
On en déduit quelques conséquence en théorie de la causalité.

\vfill
\hrule
\vs
{The research leading to these results has received funding from the European Research Council
	under the European Union's Seventh Framework Programme (FP/2007-2013) / ERC Grant
	Agreement  307062. Stefan Suhr is supported by the SFB/TRR 191 `Symplectic Structures in Geometry, Algebra and Dynamics', funded by the DFG.}

\end{small}

\newpage

\hspace{1cm}

\section{Definitions and Results}

In the setting of Lorentzian geometry, the relation between continuous time functions 
and smooth time functions, or even better temporal functions, was not well understood before the paper of Bernal and Sanchez
\cite{besa3}.
In the present paper, we study the more general question of the approximation of continuous causal 
functions by smooth ones, in the very general setting of closed cone fields, and without any causality assumption.
We introduce a class of causal functions, called  strictly causal functions, which can be approximated by smooth causal,
and even smooth Lyapunov, functions. Examples show that some causal functions do not have this property.
Yet time functions are strictly causal hence can be approximated by temporal functions, which generalizes the result of Bernal and Sanchez to the context of closed cone fields.

These considerations also have several consequences in causality theory. We recover in particular the equivalence between 
stable causality, the existence of a temporal function, and $\mK$-causality.

Although it may seem less relevant for physical applications to study noncausal cone fields, we point out that our result
also applies to the study of dynamical systems, where the cone field is reduced to a half-line. 
The relation between recurrence and Lyapounov or causal functions, respectively,  for such systems was first studied by Conley.
To our knowledge, our results are also new in this case.
%
%
%With the recent results \cite{BS} on Lyapunov functions and the stable relation of closed cone fields the question has appeared whether the stable relation 
%is the smallest closed and transitive relation  $\mK^+$ containing the causal relation of the cone field. The question was answered positively under the assumption of an 
%empty stably recurrent set in \cite{miK} for Lorentzian metrics and in \cite{mi17} for a class of closed cone fields slightly more restricted than the one considered here. A 
%counterexample was given in \cite{mi2}, see Example \ref{E1}. 
%
%The stable relation is equal to $\mK^+$ if every causal function, i.e. continuous and monotone along causal curves, 
%is approximated uniformly on compact sets by Lyapunov functions. It is known that $\mK^+$ is generated by the causal functions \cite{aus,levin}.
%In this article we will introduce a subclass of causal functions (special causal) which are approximated by Lyapunov functions. Further we discuss the consequences for causality theory of closed cone fields.
%

We work on a  Riemannian manifold $M$. A {\it convex cone} of the vector space $E$ is a convex subset $C\subset E$ such that $tx\in C$ for each $t> 0$ and $x\in C$.
The convex cone $C$ is called {\it regular} if it is not empty and  contained in an open half-space, or equivalently
 if there exists a linear form $p$ on $E$ such that $p\cdot v>0$ for each  $v\in C$.
The full cone $C=E$ will be called the {\it singular} cone.

\begin{defn}
We say that  $C\subset E$ is a \text{closed cone} if it is a convex cone which is either singular or regular and if
	$C\cup \{0\}$ is a closed subset of $E$.
	\end{defn}

Note that the empty set is a closed cone. We call it {\it degenerate}. 

A {\it cone field} $\mC$ on the manifold $M$ is a subset of the tangent bundle $TM$ such that $\mC(x):=T_xM \cap \mC$
is a convex cone for each $x$.
%We shall only use closed cone fields:

\begin{defn}
	We say that  $\mC\subset TM$ is a {closed cone field} if it is  a cone field such that the union of $\mC$  and the image of the zero section in $TM$
	is a closed subset of $TM$ and such that $\mC(x)$ is a closed cone for each $x$.
\end{defn}

A point $x\in M$ is called nondegenerate for $\mC$ if $\mC(x)\neq \emptyset$. 
The  {\it domain} of $\mC$ is the set of nondegenerate points. It is denoted by $\mD(\mC)$.
The domain of a closed  cone field is closed.
A cone field is  called {\it nondegenerate} if all points are nondegenerate, \textit{i.e.} if $\mD(\mC)=M$.
The set of singular points (the points $x$ such that $\mC(x)=T_xM$) of a closed  cone field is closed.

For a closed cone field $\mC$, recall from \cite{BS} that the curve  $\gamma:I\lto M$ is $\mC$-causal  (or just causal) if
it is locally Lipschitz and   if the inclusion $\dot \gamma(t) \in \mC(\gamma(t))\cup \{0\}$ holds for almost all $t\in I$.
The strictly causal  future  $\mJ_{\mC}^{++}(x)$ of $x$ is the set of points $y\in M$ such that there exists a 
nonconstant  causal curve  $\gamma:[0,T]\lto M$
satisfying $\gamma(0)=x$ and $\gamma(T)=y$.
The causal future  $\mJ_{\mC}^+(x)$ of $x$ is the set 
$\mJ^+_{\mC}(x)=\mJ_{\mC}^{++}(x)\cup \{x\}$.
The causal past $\mJ_{\mC}^-(x)$ is the set of points $x'\in M$ such that $x\in \mJ^{+}_{\mC}(x')$.
More generally, for each subset $A\subset M$, we denote by $\mJ^{\pm}_{\mC}(A):= \cup_{x\in A}\mJ^{\pm}_{\mC} (x)$ the causal future and past
of $A$. 
We have the inclusion $\mJ_{\mC}^+(y)\subset \mJ_{\mC}^+(x)$ if 
$y \in \mJ_{\mC}^+(x)$.

\begin{defn}
The function $f\colon M\to \RR$ is said to be causal if it is continuous, and has the property
that $f(y)\geq f(x)$ for each $y\in \mJ^+_\mC(x)$.
We denote by $\gC$ the set of causal functions.
\end{defn}

Given a causal function $f$, the point $x$ is said to be {\it neutral} if it is a singular point of $\mC$ or  if there exists
$y\in  \mJ^+_{\mC}(x)-\{x\}$ such that $f(y)=f(x)$. Otherwise, $x$ is said to be {\it strict}.

We say that $a\in \Rm$ is a {\it neutral value} of $f$ if there exists a neutral point 
$x$ such that $f(x)=a$. Otherwise, $a$ is called a {\it strict value}.

We say that $f$ is a {\it special causal function} if it is causal, and if
the set of strict values of $f$ is dense in $\Rm$. We denote by 
$\gS$ the set of special causal functions. It is classical in the Conley theory of dynamical systems to consider this kind of hypothesis. The terminology
 {\it neutral values}
is taken from \cite{FP}.

We say that $\tau$ is a Lyapunov function if it is smooth and
satisfies $d\tau _x \cdot v >0$ for each regular point $x$ of $\tau$ and 
each $v\in \mC(x)$.
We denote by $\gL$ the set of Lyapunov functions.
 We have 
$$
\gL\subset \gS\subset \gC.
$$
The first inclusion following from:
\begin{lem}
	Each Lyapunov function $\tau$ is special causal. Moreover, the set of neutral points of $\tau$ is contained in the set of critical points of $\tau$.
\end{lem}

\proof
If $\tau$ is a Lyapunov function and $\gamma$ is a causal curve, then 
the function $\tau \circ \gamma$ is locally Lipschitz and satisfies 
$(\tau\circ \gamma)'(t)\geq 0$ almost everywhere, hence it is nondecreasing.
This implies that the function $\tau$ is causal.
If there exists $t>0$ and a nonconstant causal curve $\gamma\colon [0,t]\to M$ such that $\tau\circ \gamma(t)=\tau\circ \gamma(0)$,
then $\tau\circ \gamma$ is constant on $[0,t]$, which implies that $\gamma(0)$
is a critical point of $\tau$. We deduce that the neutral points of $\tau$ are 
necessarily critical points, hence that neutral values are critical values.
By Sard's Theorem, the set of regular values of $\tau$ is dense, which implies that
the set of strict values is dense.
\qed 

%Second we have an easy consequence of the definitions.
%
%\begin{lem}\label{lem-cone}
%If $f\in \gC$ and $\tau \in \gS$ then we have $f+\tau\in \gS$.
%\end{lem}
%
%\proof
%Let $x\in M$ be a strict point for $\tau$, i.e. $x$ is not a singular point of $\mC$ and $\tau(y)>\tau(x)$ for all $y\in  \mJ_{\mC}^+(x)\setminus\{x\}$. 
%Then $x$ is strict for $f+\tau$ as well since $f(y)\ge f(x)$ for all $y\in \mJ_{\mC}^+(x)$.
%\qed
%

We define the stably recurrent set $\mR$ of $\mC$ as the intersection of the sets of critical points of Lyapunov functions.
This set is studied in \cite{BS}. In the present setting, it is proved there that:

\begin{prop}\label{prop-r}
	If $x\in \mR$, then either $x$ is a singular point of $\mC$, or there exists $y\in\mJ^+_\mC(x)-\{x\}$
	such that $\tau(y)=\tau(x)$ for each $\tau\in \gL$. As a consequence
	$$
	\mR=\bigcap _{\tau \in \gL} \textnormal{neut}(\tau)=\bigcap _{\tau \in \gL} \textnormal{crit}(\tau).
	$$
	Moreover, there exists a Lyapunov function $\tau_0$ such that 
	$$
	\mR=\textnormal{neut}(\tau_0)=\textnormal {crit}(\tau_0).
	$$
\end{prop}

We now state our main result, which will be proved in Section \ref{section2}:

\begin{thm}\label{thm1}
For each function $f\in \gS$ and each $\epsilon>0$, there
	exists a Lyapunov function $\tau\in \gL $ such that $|\tau-f|(x)<\epsilon$ for each $x$, and which is regular 
	outside of $\mR$.
\end{thm}

To our knowledge, this statement is new even in the classical Lorentzian setting
and in the setting of vector fields.
See however \cite{cgm} for results on smoothing causal functions. 

The following example, taken from \cite[Example 5.1]{mi2}, illustrates that the set $\gS$ could not be replaced by $\gC$ 
in the above statement.

\begin{exmp}\label{E1}
We consider the manifold $M=\Sm^1\times \Rm$,  and the constant cone field 
$\mC(\theta, x)=[0,\infty)\times [0,\infty)-\{(0,0)\}$. 
The function $f(\theta,x)$ on $M$ is causal if and only if it is of the 
form $f(\theta,x)=g(x)$ for some nondecreasing function $g:\Rm\lto \Rm$.
For such a causal function, all the values in the interval $g(\Rm)$ are neutral.
As a consequence, the function $f$ is special causal if and only if it is constant.
\end{exmp}

Obviously, $\gS$, and hence $\gL$, is not dense in $\gC$ in this example.

\begin{cor}
	$\mR=\bigcap _{f \in \gS}\textnormal{neut}(f)$
	\end{cor}

\proof
Since $\gL\subset  \gS$, we have the inclusion $\bigcap _{f \in \gS}\textnormal{neut}(f)\subset \mR$ by Proposition \ref{prop-r}.
Conversely, consider $x\in \mR$. 
If $x$ is a singular point of $\mC$, then it is a neutral point for each causal function.
If $x$ is not a singular point of $\mC$, there exists $y\in \mJ^+_{\mC}(x)-\{x\}$ such that $\tau(y)=\tau(x)$ for each $\tau\in \gL$.
The density of $\gL$ in $\gS$ then implies that $f(y)=f(x)$ for each special causal function $f$, hence $x$ is 
neutral for $f$.
\qed

In general the set of neutral points of a causal function can be a strict subset of $\mR$.
We recall some definitions before presenting an example.
We say that $f$ is a {\it time function} if it is a causal function without neutral points. 
Obviously, time functions are special causal.
Note that the existence of a time 
function with our definitions prevents the presence of singularities of $\mC$.
We say that $\tau$ is a {\it temporal function} if it is a Lyapunov function without critical points.

\begin{exmp}
Consider $M=S^1\times \RR^2$ with coordinates $(\theta,x,z)$. Equip $M$ with 
the cone field
$$
\mC(\theta,x,z):= \{(v_{\theta}, v_x,v_z): v_x\geq 0 \textnormal{ and } v_{\theta}v_{x}\geq z^2 v_{\theta}^2+v_z^2\}.
$$
We claim that 
$$\bigcap_{f\in\gC}\textnormal{neut}(f)=S^1\times\RR\times\{0\}
\neq \mR=M.$$
\end{exmp}
\proof
Under the embedding $(\theta,x)\mapsto (\theta,x,0)$ of $S^1\times \RR$ into $M$ the cone field $\mC$ pulls back to the cone field considered in Example
\ref{E1}.
The coordinate function $x\colon M\to\RR$ is a causal function of $(M,\mC)$, but not a temporal one.
 As seen in Example 
\ref{E1} the set of neutral points of every causal function includes $S^1\times\RR\times \{0\}$. Further the neutral points of the function  $x$
are precisely $S^1\times\RR\times \{0\}$. Therefore we have
$$\bigcap_{f\in\gC}\textnormal{neut}(f)=S^1\times\RR\times\{0\}\subset \mR.$$
In particular, all Lyapunov functions are critical, hence constant, on $S^1\times\RR\times\{0\}$.
The curves $(\theta(t),x(t),z(t))=(\theta_0+t, x_0+2t+(t_0+t)^3-t_0^3,t_0+ t)$ 
and $(\theta(t),x(t),z(t))=(\theta_0+t, x_0+2t+(t-t_0)^3+t_0^3,t_0- t)$ 
are causal.
This implies that for every $p\in M$ we have 
$$\mJ^\pm_\mC(p)\cap (S^1\times\RR\times\{0\})\neq \emptyset.$$
As a consequence,  all Lyapunov functions are constant on $M$, hence   $\mR=M$.
\qed

We recall that a time function is said to be Cauchy if it has the property that the image of $f\circ \gamma$ 
is the whole real line for each inextendible causal curve $\gamma$.

\begin{cor}\label{C1}
	There exists a time function if and only if the stably recurrent set is empty.
In this case,  each causal function $f$ can be uniformly approximated by temporal functions $\tau$.
	
	If $f$ is a Cauchy time function, then so is the temporal function $\tau$ (in this case, $\mC$ is globally hyperbolic
	in the sense of \cite{BS}).
	\end{cor}

\proof
The first claim follows directly from the definition and Proposition \ref{prop-r}.

Let us next consider a causal function $f$. Since $\mR$ is empty, there exists a time function $g'$.
The function $g:= \frac{2\epsilon}{\pi} \arctan \circ g'$ is a time  function taking values in $]-\epsilon, \epsilon[$.
The function $f+g$ is then a time function, hence it belongs to $\gS$. 
By Theorem \ref{thm1}, there exists a temporal function $\tau$ such that $|\tau-(f+g)|<\epsilon$,
hence $|\tau-f|< 2\epsilon$.

If $f$ and $g$ are two time functions such that $f-g$ is bounded, and if $f$ is Cauchy, then so is $g$.
Indeed, for each inextendible causal curve $\gamma$, we have $\sup (f\circ \gamma)=+\infty$,
hence  $\sup (g\circ \gamma)=+\infty$ and 
$\inf (f\circ \gamma)=-\infty$,
hence  $\inf (g\circ \gamma)=\infty$.
\qed

Let us discuss the implications of our result in terms of causality relations.
The stable (also called Seifert) relation $\mF^+(x)$ is defined by
$$
\mF^+(x):=\{y: \tau(y)\geq \tau(x) \;\forall \tau\in \gL\}
$$
This relation can also be defined in terms of cone enlargements, see  \cite{BS} 
for the present generality.
It is an obvious consequence of Theorem \ref{thm1} that 
$$
\mF^+(x)=\{y: \tau(y)\geq \tau(x) \;\forall \tau\in \gS\}.
$$
On the other hand, it is a general fact that 
$$
\mK^+(x)=\{y: \tau(y)\geq \tau(x) \;\forall \tau\in {\gC}\}\subset {\mF^+}(x)
$$
is the smallest transitive closed relation containing $\mJ^+_{\mC}$. 
This follows from the Auslander-Levin Theorem, \cite{aus,levin}, \cite[Theorem 2]{AA},
as was noticed in 
\cite[Theorem 3]{mi1}.

In Example \ref{E1} we have 
$\mK^+(\theta,x)=\mJ^+_{\mC}(\theta, x)= \Sm^1\times [x,\infty)$
and $\mF^+(\theta,x)=M$ for all $(\theta,x)\in M$.
This shows that equality does not hold in general between $\mK^+$ and $\mF^+$.
In the case where $\mR$ is empty, however,  we have by Corollary \ref{C1}:

\begin{prop}
	If $\mR$ is empty, $\mK=\mF$. 
\end{prop}

The equality $\mK=\mF$ is proved in \cite{miK} in the classical Lorentzian setting, and very recently \cite[Theorem 85]{mi17} 
under  assumptions slightly stronger than the present ones.

Let us finally recover and slightly generalize some results of Minguzzi on
the equivalence between $\mK$-causality and stable causality.
We say that the cone field is {\it $\mK$-causal} if $\mK^+(x)\cap \mK^-(x)=\{x\}$ for each $x$.
$\mK$-causality implies that for all $x$ and all $y\in \mJ^+_{\mC}(x)- \{x\}$ there exists $f\in \gC$
with $f(y)>f(x)$. 
Recall that the cone field is {\it stably causal} if $\mR$ is empty.
This is equivalent to {\it $\mF$-causality}, that is to the property that  
$\mF^+(x)\cap \mF^-(x)=\{x\}$ for each $x$, as follows from Proposition \ref{prop-r},
see  also \cite{Sei,HS,mi09,BS}.

\begin{prop}
	The following properties are equivalent:
	\begin{itemize}
		\item The cone field $\mC$ is stably causal.
		\item There exists a temporal function.
		\item There exists a time function.
		\item The cone field is $\mK$-causal.
	\end{itemize}
\end{prop}

\proof
We have already proved the equivalence of the three first points.

Since $\mK^+(x)\cap \mK^-(x)\subset \mF^+(x)\cap \mF^-(x)$, any of the first three conditions implies the last point.

Conversely, 
let us consider the set $\gC([0,1])$ of causal functions taking values in $[0,1]$, endowed with 
the topology of uniform convergence on compact sets. It is a separable metric space,
as a subset of the separable metric space of all continuous functions 
with the topology of uniform convergence on compact sets.
We consider a dense sequence $f_i\in \gC([0,1])$ and 
$f:=\sum_{i\geq 1} 2^{-i} f_i$, which is an element of $\gC$.
If the cone field is $\mK$-causal, then $f$ is a time function.
Indeed, if $y\in\mK^+(x)-\{x\}$, then $f_i(y)\geq f_i(x)$
for each $i$, and there exists $j$ such that $f_j(y)>f_j(x)$.
As a consequence, we get that $f(y)>f(x)$.
\qed

\section{Proof of Theorem \ref{thm1}}\label{section2}

We need more definitions from \cite{BS}.

We say that $\mE\subset TM$ is an {open cone field} if it is a  cone field  which is open as a subset of $TM$.
Then $\mE(x)$ is an open cone for each $x$. We say that $\mE$ is an open enlargement of $\mC$
if it is an open cone field containing $\mC$.

Given an open cone field $\mE$, we say that the curve $\gamma:I\lto M$ is $\mE$-timelike  (or just timelike) if
it is piecewise smooth and  if $\dot \gamma(t) \in \mE(\gamma(t))$ for  all $t$ in $I$ (at nonsmooth points, the inclusion is required to hold for 
left and right derivatives).
The  {\it chronological future}  $\mI_{\mE}^+(x)$ of $x$ is the set of points $y\in M$ such that there exists a nonconstant  timelike curve $\gamma:[0,T]\lto M$
satisfying $\gamma(0)=x$ and $\gamma(T)=y$.
The {{\it chronological past} $\mI^-_{\mE}(x)$ of $x$} is the set of points $x'\in M$ such that $x\in \mI^{+}_{\mE}(x')$.
More generally, for each subset $A\subset M$, we denote by $\mI^{\pm}_{\mE}(A):= \cup_{x\in A}\mI^{\pm}_{\mE} (x)$ the {  {\it chronological future and past
		of}} $A$. 
	The chronological future can also be defined using only smooth timelike curves, see
	\cite[Lemma 22]{BS}.

\begin{prop}\label{prop-step}
	Let $f$ be a special causal function, and let $a\in \RR$ a strict value of $f$.
	Then for each open interval $]a_-,a_+[\ni a$, there exists a Lyapunov function 
	$\tau : M\lto [a_-, a_+]$ such that 
	$\tau(x)=a_+$ if 
	$f(x)\geq a_+ $
	and $\tau(x)=a_-$ if $f(x)\leq  a_-$ and all values in $]a_-,a_+[$ are regular for $\tau$.
\end{prop}

{\sc Proof of Theorem \ref{thm1} from Proposition \ref{prop-step}:}
We fix a special causal function $f$ and $\epsilon>0$.
For each $k$, there exists a Lyapunov function $\tau_k: M\lto [k\epsilon, (k+1)\epsilon]$ such that 
$\tau_k=(k+1)\epsilon  $ on $\{f\geq (k+1)\epsilon\}$ and 
$\tau_k={k}\epsilon $ on $\{f\leq k \epsilon\}$.
The function $\tau$ which is equal to $\tau_k$ on the set 
$f^{-1}([k\epsilon, (k+1)\epsilon])$ for each $k$ is a Lyapunov function,
whose set of critical values is $\epsilon \Zm$, and which satisfies 
$|f-\tau|<\epsilon $.

In order to control the critical set, we consider a Lyapunov function $\tau_0:M\lto [0,\epsilon]$ which is
regular outside of $\mR$. Then $\tau+\tau_0$ is regular outside of $\mR$
and $|f-(\tau+\tau_0)|<2\epsilon $.
\qed

\noindent
{\sc Proof of Proposition \ref{prop-step}:}
Let $a\in \RR$ be a strict value of the special causal function $f\colon M\to \RR$. Then $U:=\{f>a\}$ is open with $\mJ^{++}_{\mC}(\overline{U})\subset U$ and 
$\partial U$ consists only of regular points of $\mC$.
Especially there exists a neighborhood of $\partial U$ consisting only of regular points of $\mC$ since the 
singular set of $\mC$ is closed. Recall from \cite{BS} the definition of a trapping domain: 

\begin{defn}
	An open set $A\subset M$ is a { \it trapping domain} for the open cone field $\mE$
	if $\mI^+_{\mE}(A)\subset A$. $A$ is a trapping domain for the closed cone field $\mC$
	if it is a  trapping domain for some open enlargement $\mE$ of $\mC$.
	We say that the domain is locally trapping at a point $x \in \partial A$ if there exists an open neighborhood
	$U$ of $x$ such that $A\cap U$ is a trapping domain for $\mC_{|U}$.
\end{defn}

We proved in \cite[Lemma 35]{BS} that a domain is trapping if and only if it is locally trapping at each point of its boundary.

\begin{prop}\label{prop-reg}
	Let $U$ be an open set such that $\mJ^{++}_{\mC}(\bar U)\subset U$ and assume that $\partial U$ consists only of regular points of $\mC$. Further let $F_i$
	be a closed set contained in $U$, and let $F_e$ be a closed set disjoint from $\bar U$.
	Then there exists a  trapping domain $A$ such that $F_i\subset A$ and $F_e$
	is disjoint from $\bar A$.
\end{prop}

\noindent
{\sc Proof of Proposition \ref{prop-reg}:}
Set $W:=M- \bar U$. For an open cone field $\mE$ denote with $\mE_W$ the open cone field which is equal to $\mE$ on $W$ 
and empty on $\bar U$.
We denote by $\partial^- \mI^+_{\mE_{W}}(y)$ the part of the boundary of $\mI^+_{\mE_{W}}(y)$ which is contained  in $W$.
Note that the open set $\mI^+_{\mE_{W}}(y)$ is locally trapping  at each point of $\partial^-\mI^+_{\mE_{W}}(y)$ and 
$\mI^+_{\mE_{W}}(y)\subset W$.

\begin{lem}\label{lem-reg}
	For each $x\in \partial U$, there exists $y_x\in W\cap \mJ^-_{\mC}(x)$ and an enlargement $\mE^x$ such that $V_x:=\mI^+_{\mE^x_{W}}(y_x)$ is bounded in $M$ and 
	disjoint from $F_e$.
\end{lem}

\noindent
{\sc Proof of Lemma \ref{lem-reg}: }
Let us fix $x\in \partial U$.
We consider a bounded open neighborhood $V$ of $x$ in $M$, which is disjoint from $F_e$ and does not contain any singular point of $\mC$.
We furthermore assume that there exists a temporal function $\tau$ on $V$, which satisfies $d\tau_z\cdot v>2|v|_z$ for each $z\in V$
and $v\in \mC(z)$.  
We will prove that there exists $y^x\in W\cap \mJ^-_{\mC}(x)$ and an open enlargement $\mE^x$ of $\mC$ such that  $\mI^+_{\mE^x_{W}}(y^x)\subset V$.
Since $V$ is bounded so is 
$\mI^+_{\mE^x_{W}}(y^x)$. Further $\mI^+_{\mE^x_{W}}(y^x)$ is disjoint from $F_e$ because $V$ is.

We consider a nested sequence of open enlargements $\mE^n$ with $\mC\subset \mE^n\subset  \mE$ and $\cap_n\mE^n=\mC$. 
The existence of such a sequence is proved in \cite[Lemma 20]{BS}. We set  $L=d(x, \partial V)/2$ where $d$ denotes the
distance with respect to the Riemannian metric.

We first prove the existence of a $\mC$-causal curve $\eta: [-\delta,0]$,
parametrized by arclength, such that $\eta(0)=x$
and $\eta(t)\in W$ for each $t\in [-\delta, 0[$.
For this, we consider a sequence $\eta_n:[-L,0]\lto V$ of $\mE_n$-timelike curves
parametrized by arclength. At the limit (possibly along a subsequence),
we obtain a curve $\tilde \eta:[-L,0]\lto V$, which is $\mC$ causal 
according to the limit curve Lemma 
\cite[Lemma 23]{BS}. 

Since $\cap \mE^n=\mC$,  the curves $\eta_n$ satisfy $d\tau_{\eta_n(t)}\cdot \dot\eta_n(t)\geq |\dot\eta_n(t)|=1$ for each $t\in [-L,0]$ when $n$ is large enough,
hence $\tau\circ \eta_n(-t)\leq \tau\circ \eta_n(0)-t$ for each $t\in [0,L]$.
At the limit, we deduce that the curve $\tilde \eta$ is not constant.
 We parametrize it by arclength and get $\eta$.
Since $\mJ^{++}_{\mC}(\bar U)\subset U$ and $\eta(0)\not \in U$,
we obtain that $\eta(t)\in W$ for $t<0$.

We now consider a sequence $y_n\in \mJ^-_{\mC}(x)\cap W$ converging to $x$.
These points can be taken on the curve $\eta$.
We claim that $\mI^+_{\mE^n_{W}}(y_n)$ is contained in $V$ when $n$ is large enough.
If we assume the contrary,
there exists  for each $n$ large enough a smooth  $\mE^n$-timelike curve  $\gamma_n\colon [0,L]\lto V\cap W$ parametrized by arclength
and starting at $y_n$.
As above, these curves converge uniformly to a nonconstant  limit 
$\gamma:[0,L]\lto  \overline W$ 
which satisfies $\gamma(0)=x$.
The fact that $\gamma(L)$ belongs to $\overline W\subset M- U$ is in contradiction with 
the hypothesis   $\mJ^{++}_{\mC}(x) \subset U$.
\qed

The set $A:=\cup_x V_x \cup \bar U$ is open in $M$. The boundary of $A$ is contained in $\cup _x \partial^- V_x$, hence $A$ is locally trapping at each point 
of its boundary, hence it is trapping by \cite[Lemma 35]{BS}. Since $\overline U\subset A$ we know that $F_i\subset A$. Further by construction $\overline A$ is disjoint from $F_e$. This ends the proof of Proposition \ref{prop-reg}.
\qed

In order to finish the proof of Proposition \ref{prop-step} we invoke two propositions from \cite{BS}: 

\begin{prop}
	Let $A$ be a trapping domain, let $F_i$ be a closed set contained in $A$, let $F_e$
	be a closed set disjoint from $\bar A$, and let $\theta$ be a point in the boundary of $A$.
	
	Then there exists a smooth (near $\mD(\mC)$) trapping domain $A'$ which contains $F_i$, whose boundary contains 
	$\theta$, and whose closure is disjoint from $F_e$.
\end{prop}

\begin{prop}
	Let $A$ be  smooth trapping domain, then 
	there exists a (smooth) Lyapunov function $\tau':M\lto [-1,1]$
	such that $A=\{\tau' >0\}$ and
	all values in $]-1,1[$ are 
	regular values of $A$
	(hence $\partial A=\{\tau'=0\}$).
	
	If $F_i$ and $F_e$ are closed sets  
	contained in $A$ and disjoint from $\bar A$, respectively, we can moreover impose that $\tau'=1$ on $F_i$ and 
	$\tau'=-1$ on $F_e$.
\end{prop}

Here we consider $F_i:=\{f\geq a_+\}$ and $F_e:=\{f\leq a_-\}$. The Lyapunov function $\tau$ 
claimed in Proposition \ref{prop-step} is defined as 
$$\tau :=\frac{1}{2}[(a_+-a_-)\tau'+(a_+ +a_-)]. $$
This ends the proof of Proposition \ref{prop-step}.

\end{document}